\pdfoutput=1

\documentclass[a4paper,12pt]{article}

\usepackage[utf8]{inputenc}
\usepackage[T1]{fontenc}
\usepackage{microtype}
\usepackage[top=2.65cm, bottom=3.7cm, left=2.7cm, right=2.7cm]{geometry}
\usepackage[charter, uppercase=upright, greeklowercase=upright]{mathdesign}
\usepackage[scaled=.96]{XCharter}
\usepackage{xparse}
\usepackage{mathtools}
\usepackage{amsthm}
\usepackage{stmaryrd}
\usepackage{hyperref}

\mathtoolsset{lastline-preskip=0pt}
\allowdisplaybreaks
\hypersetup{colorlinks, linkcolor=black, citecolor=red, urlcolor=blue}
\newcommand{\signature}[1]{\vspace{7pt}
  {\setlength{\parindent}{0pt} #1}}

\renewcommand{\epsilon}{\varepsilon}
\renewcommand{\phi}{\varphi}
\newcommand{\N}{{\mathbb{N}}}

\newcommand{\Q}{{\mathbb{Q}}}
\newcommand{\R}{{\mathbb{R}}}
\newcommand{\Ri}{\R \cup \{-\oo\}}
\newcommand{\C}{{\mathbb{C}}}
\newcommand{\D}{{\mathbb{D}}}
\newcommand{\Db}{\bar{\D}}
\newcommand{\K}{{\mathbb{K}}}
\newcommand{\Bc}{{\mathcal{B}}}
\newcommand{\Cc}{{\mathcal{C}}}
\newcommand{\Ec}{{\mathcal{E}}}
\newcommand{\Hc}{{\mathcal{H}}}
\newcommand{\Lc}{{\mathcal{L}}}
\newcommand{\Nc}{{\mathcal{N}}}
\renewcommand{\d}{{\mathrm{d}}}
\newcommand{\Sg}{[-1 ; 1]}
\newcommand{\oo}{\infty}
\newcommand{\bs}{\backslash}
\newcommand{\ra}{\rightarrow}
\newcommand{\ti}{\widetilde}
\renewcommand{\bar}{\overline}
\newcommand{\ar}{\overrightarrow}
\NewDocumentCommand{\abs}{O{} m}{#1 \lvert #2 #1 \rvert}
\NewDocumentCommand{\norm}{O{} m}{#1 \lVert #2 #1 \rVert}
\newcommand{\ndot}{\norm{\cdot}}
\newcommand{\lb}{\llbracket}
\newcommand{\rb}{\rrbracket}
\newcommand{\st}{\; | \;}
\newcommand{\tgamma}{\ti{\gamma}}
\newcommand{\trho}{\ti{\rho}}
\newcommand{\den}{\text{\upshape cnt}}
\newcommand{\Psh}{Plurisubharmonic}
\newcommand{\psh}{plurisubharmonic}
\newcommand{\PSH}{\textsc{psh}}
\newcommand{\Pshe}{Plurisubharmonicity}
\newcommand{\pshe}{plurisubharmonicity}
\newcommand{\scs}{upper semicontinuous}
\newcommand{\sce}{semicontinuity}
\newcommand{\pc}{pseudoconvex}
\newcommand{\smmp}{strong maximum modulus principle}

\DeclareMathOperator{\real}{\Re\mspace{-1mu}\mathfrak{e}}
\DeclareMathOperator*{\essup}{ess \; sup}

\theoremstyle{definition}
\newtheorem*{dfn}{Definition}
\newtheorem*{ntn}{Notation}

\theoremstyle{plain}
\newtheorem{thm}{Theorem}
\newtheorem{prp}[thm]{Proposition}
\newtheorem{lem}[thm]{Lemma}
\newtheorem{cor}[thm]{Corollary}

\theoremstyle{remark}
\newtheorem*{rmq}{Remark}
\newtheorem*{rms}{Remarks}
\newtheorem*{exm}{Example}

\title{Convexity, \pshe\ and the \smmp\ in Banach spaces}
\author{Anne-Edgar \textsc{Wilke}}
\date{}

\begin{document}

\maketitle

\textsc{Abstract.}~-- In this article, we first try to make the known
analogy between convexity and \pshe\ more precise. Then we introduce a
notion of strict \pshe\ analogous to strict convexity, and we show how
this notion can be used to study the \smmp\ in Banach spaces. As an
application, we define a notion of $L^p$ direct integral of a family of
Banach spaces, which includes at once Bochner $L^p$ spaces, $\ell^p$
direct sums and Hilbert direct integrals, and we show that under
suitable hypotheses, when $p < \oo$, an $L^p$ direct integral
satisfies the \smmp\ if and only if almost all members of the family
do. This statement can be considered as a rewording of several known
results, but the notion of strict plurisubharmonicity yields a new
proof of it, which has the advantage of being short, enlightening and
unified.

\section{Introduction}

\subsection{Convexity and \pshe}

\Psh\ functions (\PSH, for short) were introduced independently by Oka
\cite[p.~40]{O42} and Lelong \cite[p.~306, déf.~1]{L45}. Since the
origins of the theory, it has been observed that there is a certain
analogy between convex functions and \PSH\ functions; in fact, Oka
called the latter \emph{\pc\ functions}. To make this analogy
apparent, Bremermann \cite[pp.~34-38]{B56} collected a list of
properties satisfied by convex functions, and showed that for each of
them, there is a corresponding property satisfied by \PSH\ functions.

Here are Bremermann's main ideas, slightly reformulated. A continuous
function $f : \R \ra \R$ is said to be convex, or \emph{sublinear}, if
for all compact intervals $I \subset \R$ and for all affine functions
$\alpha : I \ra \R$, the inequality $f \leq \alpha$ holds on $I$ as
soon as it holds on the boundary of $I$. A continuous function $f :
\R^n \ra \R$ is said to be convex, or \emph{plurisublinear}, if for
all affine maps $\gamma : \R \ra \R^n$, the composition $f \circ
\gamma$ is convex.

In the same way, an \scs\ function $f : \C \ra \Ri$ is said to be
subharmonic if for all connected, smoothly bounded compact sets $K
\subset \C$ and for all functions $\alpha : K \ra \R$ continuous and
harmonic in the interior of $K$, the inequality $f \leq \alpha$ holds
on $K$ as soon as it holds on the boundary of $K$. An \scs\ function
$f : \C^n \ra \Ri$ is said to be \psh\ if for all affine maps $\gamma
: \C \ra \C^n$, the composition $f \circ \gamma$ is subharmonic.

As Bremermann remarks, convexity and \pshe\ can be defined more
quickly in the following way. A continuous function $f : \R^n \ra \R$
is convex if and only if for all affine maps $\gamma : \R \ra \R^n$,
\begin{equation} \label{eq:conv_intro}
  f(\gamma(0)) \leq \frac{f(\gamma(-1)) + f(\gamma(1))}{2}.
\end{equation}
This amounts to asking that the value of $f \circ \gamma$ at the
centre of the unit ball of $\R$ be less than its mean on the sphere.
Similarly, an \scs\ function $f : \C^n \ra \Ri$ is \PSH\ if and only
if for all affine maps $\gamma : \C \ra \C^n$,
\begin{equation} \label{eq:psh_intro}
  f(\gamma(0)) \leq \frac{1}{2 \pi} \int_0^{2 \pi} f(\gamma(e^{i t}))
  \, \d t,
\end{equation}
which amounts to asking that the value of $f \circ \gamma$ at the
centre of the unit ball of $\C$ be less than its mean on the sphere.

From the above discussion, it is tempting to conclude, as does
Bremermann, that the analogy between convex functions and
\PSH\ functions is obtained by replacing $\R^n$ with $\C^n$, real
affine maps with complex affine maps, and sublinearity conditions with
subharmonicity conditions. These ideas permeate much of the literature
on the subject; see for instance \cite[p.~225]{H94}. Yet we will show
that this dictionary misses an essential point and therefore is
unsatisfactory.

Indeed, a fundamental result, due to Lelong \cite[p.~325, n°~17]{L45},
states that \PSH\ functions are stable under composition with a
holomorphic map, from which one can define the notion of
\PSH\ function on a holomorphic manifold. This result does not appear
in Bremermann's list, which is understandable, since from his point of
view, there is no analogous result for convex functions.

In fact, the natural domain of a convex function is a real affine
space, while the natural domain of a \PSH\ function is a holomorphic
manifold, or even a complex analytic space\footnote{The approach taken
in this article would allow one to define the notion of convex
function more generally on any topological space $X$, as soon as one
chooses a class of continuous maps $\gamma : \Sg \ra X$ playing the
role of affine maps. An important case is when $X$ is a Riemannian
manifold and the maps $\gamma$ are geodesic segments. In the same way,
one could define the notion of \PSH\ function on any topological space
$X$ equipped with a class of maps $\gamma : \Db \ra X$ playing the
role of holomorphic maps.}. Therefore, real affine maps do not
correspond to complex affine maps, but rather to holomorphic maps.

Thus it is preferable to define the notion of \PSH\ function in the
following way: if $X$ is a holomorphic manifold, an \scs\ function $f
: X \ra \Ri$ is said to be \PSH\ if the
inequality~\eqref{eq:psh_intro} holds for all \emph{holomorphic} maps
$\gamma : \Db \ra X$, where $\Db \subset \C$ is the closed unit disc,
with the understanding that a map defined on $\Db$ is holomorphic if
it extends to a holomorphic map on a neighbourhood of $\Db$.

In the case $X = \C^n$, if $f$ satisfies this inequality for all
\emph{affine} maps $\gamma : \Db \ra \C^n$, then $f$ is \PSH: this is
the contents of Lelong's result. But this fact is best seen, not as a
definition, but rather as a characterisation, valid in a special case,
and which is far from obvious. We will not use it in this article.

Beyond the aesthetic aspect, a good understanding of the analogy
between convexity and \pshe\ enables one to obtain certain non-trivial
results about \PSH\ functions by adapting the proofs of the
corresponding, usually easier, statements concerning convex functions.
We hope that the results presented in this article will serve to
illustrate this phenomenon.

\subsection{Strict \pshe}

A continuous function $f : \R^n \ra \R$ is strictly convex if the
inequality~\eqref{eq:conv_intro} holds strictly for all non-constant
affine maps $\gamma : \R \ra \R^n$. According to Bremermann's
dictionary, this suggests the following definition: an \scs\ function
$f : \C^n \ra \Ri$ is strictly \PSH\ if the
inequality~\eqref{eq:psh_intro} holds strictly for all non-constant
affine maps $\gamma : \C \ra \C^n$. This is, with a different
formulation, Carmignani's definition \cite[pp.~285-286, def.~1.1 and
  1.2]{C76}\footnote{After changing Carmignani's definition 1.1 so
that strictly subharmonic functions are assumed to be finite on a
dense subset, and correcting definition 1.2, which erroneously omits
the hypothesis $w \neq 0$.}. However, this definition is very
unsatisfactory: indeed, we will see an example showing that the class
of functions thus obtained is not stable under composition with a
biholomorphism.

It is therefore preferable, according to the principles explained
above, to say that an \scs\ function $f : X \ra \Ri$ on a holomorphic
manifold $X$ is strictly \PSH\ if the inequality~\eqref{eq:psh_intro}
holds strictly for all non-constant \emph{holomorphic} maps $\gamma :
\Db \ra X$. We will see several results showing that strict \pshe,
thus defined, is a natural notion, analogous to strict convexity.

In the real case, there exists a stronger notion than strict
convexity: a function $f : \R^n \ra \R$ is said to be strongly convex
if it can be written locally as the sum of a convex function and a
$\Cc^2$ function $\epsilon$ such that $d^2 \epsilon$ is a positive
definite symmetric bilinear form at every point.

The analogous notion in the complex case is the following: a function
$f : X \ra \Ri$ is said to be strongly \PSH\ if it can be written
locally as the sum of a \PSH\ function and a $\Cc^2$ function
$\epsilon$ such that $\partial \bar{\partial} \epsilon$ is a positive
definite hermitian form at every point.

Unfortunately, it is a common practice in the literature to call
strongly \PSH\ functions \emph{strictly \PSH}. This situation is
unhappy, because strong \pshe\ is analogous to strong convexity, and
not to strict convexity.

\subsection{The \smmp\ in Banach spaces}

An $\R$-Banach space $E$ is said to be strictly convex if every affine
map $\gamma : \Sg \ra E$ whose image is contained in the unit sphere
is constant.

By analogy, it might be tempting to say that a $\C$-Banach space $E$
is strictly convex in the complex sense if every affine map $\gamma :
\Db \ra E$ whose image is contained in the unit sphere is constant.
This definition was strongly suggested by Thorp and Whitley
\cite{TW67}, and explicitly given by Globevnik \cite[p.~175,
  def.~1]{G75}.

However, the analogy turns out to be more satisfying if one asks
instead that every \emph{holomorphic} map $\gamma : \Db \ra E$ whose
image is contained in the unit sphere be constant. In this case, in
order to keep the terminology consistent, we will say that $E$ is
strictly \PSH.

The main result of Thorp and Whitley \cite[p.~641, th.~3.1]{TW67},
slightly reformulated, is that both definitions are actually
equivalent, that is, that $E$ is strictly \PSH\ if and only if every
\emph{affine} map $\gamma : \Db \ra E$ whose image is contained in the
unit sphere is constant. But this fact is best seen, not as a
definition, but rather as a non-trivial characterisation. We will not
use it in this article.

Beyond the analogy with strictly convex spaces, the importance of
strictly \PSH\ spaces lies in the fact that a $\C$-Banach space is
strictly \PSH\ if and only if it satisfies the \smmp, that is, if and
only if every holomorphic map from a connected manifold $X$ to $E$
whose norm has a local maximum is constant.

Strict convexity of an $\R$-Banach space can be characterised in the
following way: $(E, \ndot)$ is strictly convex if and only if for
every (or for one) increasing, strictly convex map $\psi : \R_+ \ra
\R$, the composition $\psi \circ \ndot$ is strictly convex. The
analogous statement for a $\C$-Banach space is the following: $(E,
\ndot)$ is strictly \PSH\ if and only if for every (or for one)
strictly convex map $\psi : \Ri \ra \Ri$, the composition $\psi \circ
\log \ndot$ is strictly \PSH.

These two results will give us simple characterisations of strict
convexity and strict \pshe\ of $L^p$ direct integrals, that we will
now present.

\subsection{Direct integrals}

The notion of direct integral used in this article is rather basic,
but sufficient to include at once Bochner $L^p$ spaces, $\ell^p$
direct sums and Hilbert direct integrals. One may consult
\cite[pp.~61-62]{HLR91} and \cite[pp.~683-686]{JR17} for a more
elaborate theory.

Let $(S, \Sigma, \mu)$ be a measure space, let $\Ec = (E_s)_{s \in S}$
be a \emph{measurable} family, in a sense that we will define, of real
or complex Banach spaces, and let $p \in [1 ; \oo]$. A section of
$\Ec$ is an element of the product $\prod_{s \in S} E_s$. Given a
section $\sigma$ satisfying an appropriate measurability condition,
let $\norm{\sigma}_p$ be the $p$-norm of the function $s \mapsto
\norm{\sigma(s)}_{E_s}$; explicitly,
\begin{equation} \label{eq:norm_p}
  \norm{\sigma}_p = \left\{
  \begin{aligned}
    & \left( \int_S \norm{\sigma(s)}_{E_s}^p
    \d \mu(s) \right)^\frac{1}{p}
    & \quad & \text{if $p < \oo$}, \\
    & \essup_{s \in S} \, \norm{\sigma(s)}_{E_s}
    & \quad & \text{if $p = \oo$}.
  \end{aligned}
  \right.
\end{equation}
Then $\sigma$ is said to be $p$-integrable if $\norm{\sigma}_p < \oo$,
and the $L^p$ direct integral of the family $\Ec$ is defined to be the
space of $p$-integrable sections, up to equality almost everywhere,
equipped with the norm $\ndot_p$. It is a Banach space, denoted by
$L^p(\Ec)$.

If the family $\Ec$ is constant, equal to a Banach space $E$, then
$L^p(\Ec)$ is the Bochner space $L^p(S, \Sigma, \mu ; E)$. In the case
where $E$ has dimension $1$, one recovers Lebesgue $L^p$ spaces.

Suppose that the $\sigma$-algebra $\Sigma$ is discrete, that $\mu$ is
the counting measure and that the $E_s$ are pairwise distinct. Then
$L^p(\Ec)$ is essentially the $\ell^p$ direct sum of the family $\Ec$,
denoted by $\ell^p(\Ec)$. More precisely, $L^p(\Ec)$ is the closed
subspace of $\ell^p(\Ec)$ whose elements are the sections with
countable support; this subspace coincides with $\ell^p(\Ec)$ except
when $p = \oo$ and the set of those $s \in S$ such that $E_s$ is
non-zero is uncountable.

Finally, Hilbert direct integrals correspond to the case where $p = 2$
and each $E_s$ is equal to one of the spaces $\ell^2_n$, for $n \in
\N$, or to $\ell^2_\oo$.

Here are now the results promised. For conciseness purposes, the
statements are slightly less general than what will be proved in the
article.

\begin{thm} \label{thm:conv_int_intro}
  Suppose that $\mu$ is $\sigma$-finite, that $\Ec$ is a discrete
  measurable family of $\R$-Banach spaces, and that $1 < p < \oo$.
  The direct integral $L^p(\Ec)$ is strictly convex if and only if
  $E_s$ is strictly convex for almost all $s$.
\end{thm}

\begin{thm} \label{thm:psh_int_intro}
  Suppose that $\mu$ is $\sigma$-finite, that $\Ec$ is a discrete
  measurable family of $\C$-Banach spaces, and that $1 \leq p < \oo$.
  The direct integral $L^p(\Ec)$ is strictly \PSH\ if and only if
  $E_s$ is strictly \PSH\ for almost all $s$.
\end{thm}

The notion of \emph{discrete} family which appears in these statements
is, as we will see, insignificant in practice: indeed, when $\mu$ is
$\sigma$-finite, the existence of a non-discrete measurable family
cannot be proved in ZFC.

Even if we will prove Theorems \ref{thm:conv_int_intro} and
\ref{thm:psh_int_intro} through entirely parallel methods, it is
important to note that the statements themselves are not rigorously
analogous. Indeed, definition~\eqref{eq:norm_p} is problematic in this
respect, because the true complex analogue of $\ndot_{E_s}$ is $\log
\ndot_{E_s}$, and not $\ndot_{E_s}$. An examination of the proofs
reveals that this discrepancy is the origin of the difference between
the hypothesis $1 < p < \oo$ in the real case and the hypothesis $1
\leq p < \oo$ in the complex case.

Let us now give a few immediate consequences of Theorems
\ref{thm:conv_int_intro} and \ref{thm:psh_int_intro}.

\begin{cor} \label{cor:conv_L^p}
  Suppose that $\mu$ is non-zero and $\sigma$-finite and that $1 < p <
  \oo$, and let $E$ be an $\R$-Banach space. The Bochner space $L^p(S,
  \Sigma, \mu ; E)$ is strictly convex if and only if $E$ is. In
  particular, the Lebesgue space $L^p(S, \Sigma, \mu ; \R)$ is
  strictly convex.
\end{cor}

\begin{cor} \label{cor:psh_L^p}
  Suppose that $\mu$ is non-zero and $\sigma$-finite and that $1 \leq
  p < \oo$, and let $E$ be a $\C$-Banach space. The Bochner space
  $L^p(S, \Sigma, \mu ; E)$ is strictly \PSH\ if and only if $E$ is.
  In particular, the Lebesgue space $L^p(S, \Sigma, \mu ; \C)$ is
  strictly \PSH.
\end{cor}

\begin{cor} \label{cor:conv_l^p}
  Suppose that $S$ is countable, that $1 < p < \oo$, and that the
  $E_s$ are $\R$-Banach spaces. Then $\ell^p(\Ec)$ is strictly convex
  if and only if $E_s$ is strictly convex for all $s$.
\end{cor}

\begin{cor} \label{cor:psh_l^p}
  Suppose that $S$ is countable, that $1 \leq p < \oo$, and that the
  $E_s$ are $\C$-Banach spaces. Then $\ell^p(\Ec)$ is strictly
  \PSH\ if and only if $E_s$ is strictly \PSH\ for all $s$.
\end{cor}

It is not difficult to see that Corollaries \ref{cor:conv_l^p} and
\ref{cor:psh_l^p} imply the same statements without the countability
hypothesis on $S$. Taking this remark into account, it turns out that
Theorem~\ref{thm:conv_int_intro} is a consequence of Corollaries
\ref{cor:conv_L^p} and \ref{cor:conv_l^p}, and that
Theorem~\ref{thm:psh_int_intro} is a consequence of Corollaries
\ref{cor:psh_L^p} and \ref{cor:psh_l^p}: indeed, we will see that
$L^p$ direct integrals are in fact $\ell^p$ direct sums of Bochner
$L^p$ spaces. Thus one sees that the notion of direct integral is a
means to state and prove results about $\ell^p$ direct sums and
Bochner spaces in a unified way.

Corollary~\ref{cor:conv_l^p} was proved by Day \cite[p.~314]{D41}
\cite[p.~520, th.~6]{D55} through a simple and direct method, which
can also give Corollary~\ref{cor:conv_L^p}, as the author remarks
\cite[p.~521]{D55}. This method can actually be used to prove Theorems
\ref{thm:conv_int_intro} and \ref{thm:psh_int_intro} when $1 < p <
\oo$, but probably not in the general case, as we will see in the
appendix.

Corollary~\ref{cor:psh_l^p} was proved by Jamison, Loomis and Rousseau
\cite[p.~205, th.~3.15 and p.~208, th.~4.2]{JLR85}. The case $1 < p <
\oo$ is dealt with through a method due to V. and
I. Istr\u{a}\c{t}escu \cite[p.~424, th.~2.4]{II79}, using Thorp and
Whitley's results and Lumer's theory of semi-scalar products. The
method used for $p = 1$ in fact works for all values of $p$, but on
the other hand it does not yield Corollary~\ref{cor:psh_L^p}, because
it relies on an explicit description of the dual space of
$\ell^1(\Ec)$ which does not hold in general for Bochner spaces
\cite[p.~98, th.~1]{DU77}.

Finally, Corollary~\ref{cor:psh_L^p} was proved by Thorp and Whitley
when $E = \C$, and by Dilworth \cite[p.~499, th.~2.5]{D86} in the
general case, using Thorp and Whitley's results. It is essentially
Dilworth's method that we will follow to prove
Theorem~\ref{thm:psh_int_intro}, but in a conceptual framework which
makes the proof more enlightening, shorter and less technical, and
avoids appealing to Thorp and Whitley's work.

\subsection{Organisation of the article and generalities}

Subsections 2.1 and 2.2 give a few basic facts about convex and
\PSH\ functions, emphasising the striking parallelism between the two
theories, which is apparent in definitions as well as in statements
and proofs.

Subsection 2.3 proves a form of Jensen's inequality which is used in
the next subsection. This result is known, at least in similar
contexts.

Subsection 2.4, essentially independent from the rest of the article,
contains two results which are intended to show that the notion of
strict \pshe\ is natural and analogous to strict convexity.

Section 3 deals with strictly convex and strictly \PSH\ spaces, again
insisting on the parallelism between both theories, and shows the
connection to the \smmp.

Section 4 develops the notion of direct integral, and, as an
application of the concepts introduced in the article, proves Theorems
\ref{thm:conv_int_intro} and \ref{thm:psh_int_intro}. Finally, the
appendix proves the same theorems again in the case $1 < p < \oo$
using Day's method.

Recall that a map defined on the closed unit disc $\Db \subset \C$ is
said to be holomorphic if it extends to a holomorphic map on a
neighbourhood of $\Db$.

We will use the fact that the origin of a topological vector space on
a non-discretely valued field has a fundamental system of
neighbourhoods whose members are balanced, that is, stable under
multiplication by scalars of absolute value less than $1$ \cite[I,
  p.~7, prop.~4]{EVT}.

A topological affine space is an affine space $E$ whose direction
$\ar{E}$ is a topological vector space. In this case, there is a
unique topology on $E$ such that for all $x \in E$, the bijection
$\ar{E} \ra E$ which maps $v$ to $x + v$ is a homeomorphism.

We will call a topological vector space whose topology can be defined
by a norm a \emph{normable space}.

The holomorphic manifolds considered in this article are sets equipped
with a holomorphic atlas whose charts have target an open set in a
complete $\C$-normable space. Please refer to \cite{VAR_1-7} or
\cite[pp.~7-16]{Douady} about these manifolds.

For Bochner integration theory and Bochner $L^p$ spaces, one may
consult the introductory text \cite[pp.~397-404]{Cohn} and the more
detailed expositions \cite[pp.~41-52]{DU77} and
\cite[pp.~1-30]{HNVW16}. Let us remark that Bochner theory, in full
generality, enables one to integrate functions defined on a measure
space with values in a complete \emph{normable} space.

\vspace{0.7cm}

\textsc{Acknowledgements.}~-- The author is very grateful to Karim
Belabas for his guidance throughout the preparation of this work, and
to INRIA for its financial support. Thanks also to Sébastien Boucksom
for his comments on an early version of this article, and to the
referee for his suggestions.

\section{Functions}

\subsection{Convex functions}

\begin{dfn}
  Let $E$ be a topological $\R$-affine space and $X \subset E$ a
  convex subset. A continuous function $f : X \ra \R$ is said to be
  convex if for all affine maps $\gamma : \Sg \ra E$ with values in
  $X$,
  \begin{equation} \label{eq:conv}
    f(\gamma(0)) \leq \frac{f(\gamma(-1)) + f(\gamma(1))}{2}.
  \end{equation}
  Moreover, $f$ is said to be strictly convex if the
  inequality~\eqref{eq:conv} is strict as soon as $\gamma$ is
  non-constant.
\end{dfn}

\begin{prp} \label{prp:conv_comp}
  Let $E_1$ and $E_2$ be two topological $\R$-affine spaces, $X_1
  \subset E_1$ and $X_2 \subset E_2$ convex subsets, $\phi : E_1 \ra
  E_2$ a continuous affine map such that $\phi(X_1) \subset X_2$, and
  $f : X_2 \ra \R$ a convex function. The composition $f \circ \phi$
  is convex. If $\phi$ is injective and $f$ is strictly convex, then
  $f \circ \phi$ is strictly convex.
\end{prp}

\begin{proof}
  Immediate.
\end{proof}

\begin{prp}[Maximum principle] \label{prp:conv_max}
  Let $E$ be a topological $\R$-affine space and $X \subset E$ a
  convex open subset. Every convex function $f : X \ra \R$ having a
  global maximum is constant.
\end{prp}

\begin{proof}
  Indeed, let $M$ be the maximum of $f$, and let $U = \{x \in X \st
  f(x) = M\}$. It is a non-empty closed set; if we show that it is
  also open, then we will be able to conclude that $U = X$ by
  connectedness.

  Thus, let $x \in U$, and let $V$ be a balanced neighbourhood of the
  origin in $\ar{E}$ such that $x + V$ is contained in $X$. Let $y \in
  x + V$, and let $\gamma : \Sg \ra E$ be given by $\gamma(t) = x + t
  (y - x)$; it is an affine map with values in $X$ such that
  $\gamma(0) = x$ and $\gamma(1) = y$. Since $f$ is convex,
  \begin{equation*}
    f(\gamma(0)) \leq \frac{f(\gamma(-1)) + f(\gamma(1))}{2}.
  \end{equation*}
  The left side of this inequality is $M$, and the right side is at
  most $M$, since $f(\gamma(-1))$ and $f(\gamma(1))$ are at most $M$.
  Hence we have equality; therefore, $f(\gamma(-1)) = f(\gamma(1)) =
  M$, so $f(y) = M$, which means that $y \in U$, and finally $x + V
  \subset U$.
\end{proof}

\subsection{\Psh\ functions}

\begin{dfn}
  Let $X$ be a holomorphic manifold. An \scs\ function $f : X \ra \Ri$
  is said to be \PSH\ if for all holomorphic maps $\gamma : \Db \ra
  X$,
  \begin{equation} \label{eq:psh}
    f(\gamma(0)) \leq \frac{1}{2 \pi} \int_0^{2 \pi} f(\gamma(e^{i
      t})) \, \d t.
  \end{equation}
  Moreover, $f$ is said to be strictly \PSH\ if the
  inequality~\eqref{eq:psh} is strict as soon as $\gamma$ is
  non-constant.
\end{dfn}

\begin{rmq}
  The integral is a well-defined element of $\Ri$, because $f$ is
  bounded above on every compact set by \sce.
\end{rmq}

\begin{rmq}
  In the case where $X$ is an open subset of a complete $\C$-normable
  space $E$, if $f$ satisfies the inequality~\eqref{eq:psh} for all
  \emph{affine} maps $\gamma : \Db \ra E$ with values in $X$, then $f$
  is \PSH: see \cite[p.~325, n°~17]{L45} when $E$ has finite
  dimension, and \cite[p.~172, th.~4.3]{L68} in the general case. On
  the other hand, it can happen that the inequality~\eqref{eq:psh} is
  strict for all non-constant affine maps without $f$ being strictly
  \PSH, as the following example shows.
\end{rmq}

\begin{exm}
  Let $f : \C \ra \Ri$ be a strictly \PSH\ function, and let $\pi :
  \C^2 \ra \C$ be the projection on the second factor. The composition
  $f \circ \pi$ is \PSH\ by Proposition~\ref{prp:psh_comp}, but not
  strictly \PSH, because its composition with the affine map $\gamma :
  z \mapsto (z, 0)$ is constant. Let $\phi : \C^2 \ra \C^2$ be the
  biholomorphism given by $\phi(z_1, z_2) = (z_1, z_1^2 + z_2)$. The
  preceding discussion shows that the composition $f \circ \pi \circ
  \phi$ is \PSH, but not strictly \PSH. However, if $\gamma : z
  \mapsto (a z + b, c z + d)$ is an affine map such that
  \begin{equation*}
    (f \circ \pi \circ \phi)(\gamma(0)) = \frac{1}{2 \pi}
    \int_0^{2 \pi} (f \circ \pi \circ \phi)(\gamma(e^{i t})) \, \d t,
  \end{equation*}
  then $\pi \circ \phi \circ \gamma$ is constant, since $f$ is
  strictly \PSH. Thus the map $z \mapsto (a z + b)^2 + c z + d$ is
  constant. This implies that $a = c = 0$, so $\gamma$ is constant.
\end{exm}

\begin{prp} \label{prp:psh_comp}
  Let $X_1$ and $X_2$ be two holomorphic manifolds, $\phi : X_1 \ra
  X_2$ a holomorphic map and $f : X_2 \ra \Ri$ a \PSH\ function. The
  composition $f \circ \phi$ is \PSH. If the differential of $\phi$ is
  injective outside a discrete subset of $X_1$ and $f$ is strictly
  \PSH, then $f \circ \phi$ is strictly \PSH.
\end{prp}

\begin{proof}
  The first claim is immediate. Thus, suppose that $\d \phi$ is
  injective outside a discrete set $D \subset X_1$ and that $f$ is
  strictly \PSH, and let $\gamma : \Db \ra X_1$ be a holomorphic map
  such that
  \begin{equation*}
    (f \circ \phi)(\gamma(0)) = \frac{1}{2 \pi} \int_0^{2 \pi}
    (f \circ \phi)(\gamma(e^{i t})) \, \d t.
  \end{equation*}
  Since $f$ is strictly \PSH, $\phi \circ \gamma$ is constant, so for
  all $z \in \D$, $(\phi \circ \gamma)'(z) = \d \phi(z)(\gamma'(z)) =
  0$. Thus the restriction of $\gamma$ to the open set $U = \{z \in \D
  \st \gamma'(z) \neq 0\}$ has values in $D$, which is discrete;
  therefore, $\gamma$ is constant on each connected component of $U$,
  hence its derivative is zero, so $U = \varnothing$ and $\gamma$ is
  constant on $\D$, and thus on $\Db$ by continuity.
\end{proof}

\begin{prp}[Maximum principle] \label{prp:psh_max}
  Let $X$ be a connected holomorphic manifold. Every \PSH\ function $f
  : X \ra \Ri$ having a global maximum is constant.
\end{prp}

\begin{proof}
  Indeed, let $M$ be the maximum of $f$, and let $U = \{x \in X \st
  f(x) = M\}$. It is a non-empty set, which is closed since $f$ is
  \scs\ and the condition $f(x) = M$ is equivalent to $f(x) \geq M$.
  We will see that $U$ is open, from which we will be able to conclude
  that $U = X$ by connectedness.

  Thus, let $x \in U$, and let $V$ be a balanced open neighbourhood of
  the origin in a complete $\C$-normable space and $\phi$ a
  biholomorphism between $V$ and an open neighbourhood of $x$,
  satisfying $\phi(0) = x$. Let $y \in \phi(V)$, and let $\gamma : \Db
  \ra X$ be given by $\gamma(z) = \phi(z \cdot \phi^{-1}(y))$; it is a
  holomorphic map such that $\gamma(0) = x$ and $\gamma(1) = y$. Since
  $f$ is \PSH,
  \begin{equation*}
    f(\gamma(0)) \leq \frac{1}{2 \pi} \int_0^{2 \pi} f(\gamma(e^{i
      t})) \, \d t.
  \end{equation*}
  The left side of this inequality is $M$, and the right side is at
  most $M$, since all the values $f(\gamma(e^{i t}))$ are at most $M$.
  Hence we have equality; therefore, for almost all $t \in [0 ; 2
    \pi]$, $f(\gamma(e^{i t})) = M$. By \sce, this relation actually
  holds for all $t$, so in particular $f(y) = M$, which means that $y
  \in U$, and finally $\phi(V) \subset U$.
\end{proof}

\subsection{Jensen's inequality}

\begin{lem} \label{lem:conv_aff}
  Let $E$ be a topological $\R$-affine space and $X \subset E$ a
  convex open subset, and let $f : X \ra \R$ be a convex function and
  $x_0 \in X$. There exists a continuous affine map $\alpha : E \ra
  \R$ such that $\alpha \leq f$ and $\alpha(x_0) = f(x_0)$. If $f$ is
  strictly convex, then $\alpha$ and $f$ only coincide at $x_0$.
\end{lem}

\begin{proof}
  Define
  \begin{equation*}
    C = \{(x, t) \in E \times \R \st x \in X \text{ and } f(x) < t\}.
  \end{equation*}
  It is a convex open subset of $E \times \R$, which does not contain
  the point $(x_0, f(x_0))$. According to the Hahn-Banach theorem
  \cite[II, p.~39, th.~1]{EVT}, there exists a closed hyperplane
  containing $(x_0, f(x_0))$ and not intersecting $C$. This amounts to
  saying that there exists a continuous affine map $\alpha : E \ra \R$
  and a real number $\lambda$ such that $\alpha(x_0) = \lambda f(x_0)$
  and for all $(x, t) \in C$, $\alpha(x) < \lambda t$.

  We observe that $\lambda \neq 0$: indeed, there exists $t \in \R$
  such that $(x_0, t) \in C$, so if $\lambda$ were zero, we would have
  both $\alpha(x_0) = 0$ and $\alpha(x_0) < 0$. Thus, by replacing
  $\alpha$ by $\lambda^{-1} \alpha$, one reduces to $\lambda = 1$, in
  which case the two desired properties $\alpha \leq f$ and
  $\alpha(x_0) = f(x_0)$ are satisfied.

  Let $x \in X$ such that $\alpha(x) = f(x)$, and let $\gamma : \Sg
  \ra E$ be the affine map such that $\gamma(-1) = x_0$ and $\gamma(1)
  = x$. Its image is contained in $X$. Since $\alpha$ is affine and
  $f$ is convex,
  \begin{align} \label{eq:conv_aff}
    \alpha(\gamma(0)) & = \frac{\alpha(\gamma(-1)) +
      \alpha(\gamma(1))}{2}, & f(\gamma(0)) & \leq \frac{f(\gamma(-1))
      + f(\gamma(1))}{2}.
  \end{align}
  The inequality $\alpha(\gamma(0)) \leq f(\gamma(0))$ and the
  equalities $\alpha(x_0) = f(x_0)$ and $\alpha(x) = f(x)$ show that
  the inequality~\eqref{eq:conv_aff} is an equality. Thus if $f$ is
  strictly convex, then $\gamma$ is constant, so $x = x_0$.
\end{proof}

\begin{lem}[Jensen's inequality] \label{lem:jensen}
  Let $E$ be a complete $\R$-normable space, $X \subset E$ a convex
  open subset, and $f : X \ra \R$ a convex function. Let also $(S,
  \Sigma, \mu)$ be a measure space of total mass $1$ and $\eta : S \ra
  E$ an integrable map with values in $X$. The integral $\int_S \eta
  \, \d \mu$ belongs to $X$, and
  \begin{equation*}
    f \left( \int_S \eta \, \d \mu \right)
    \leq \int_S (f \circ \eta) \, \d \mu.
  \end{equation*}
  If equality holds and $f$ is strictly convex, then $\eta$ is
  essentially constant.
\end{lem}

\begin{proof}
  Let $m = \int_S \eta \, \d \mu$. That $m$ belongs to $X$ can be seen
  by writing $X$ as an intersection of open half-spaces thanks to the
  Hahn-Banach theorem \cite[II, p.~39, th.~1]{EVT}. According to
  Lemma~\ref{lem:conv_aff}, there exists a continuous affine map
  $\alpha : E \ra \R$ such that $\alpha \leq f$ and $\alpha(m) =
  f(m)$. We have
  \begin{equation*}
    \begin{split}
      f \left( \int_S \eta \, \d \mu \right)
      & = \alpha \left( \int_S \eta \, \d \mu \right) \\
      & = \int_S (\alpha \circ \eta) \, \d \mu \\
      & \leq \int_S (f \circ \eta) \, \d \mu.
    \end{split}
  \end{equation*}
  If equality holds, then $\alpha \circ \eta$ and $f \circ \eta$
  coincide almost everywhere. In the case where $f$ is strictly
  convex, Lemma~\ref{lem:conv_aff} then shows that $\eta$ equals $m$
  almost everywhere.
\end{proof}

\subsection{Relations between convexity and \pshe}

\begin{thm} \label{thm:conv_psh}
  Let $E$ be a complete $\C$-normable space and $X \subset E$ a convex
  open subset, and let $f : X \ra \R$ be a continuous function. If $f$
  is convex, then $f$ is \PSH, and if $f$ is strictly convex, then $f$
  is strictly \PSH.
\end{thm}

\begin{proof}
  Suppose that $f$ is convex and let $\gamma : \Db \ra X$ be a
  holomorphic map. According to Cauchy's formula and Jensen's
  inequality,
  \begin{equation*}
    \begin{split}
      f(\gamma(0))
      & = f \left( \frac{1}{2 \pi} \int_0^{2 \pi} \gamma(e^{i t})
      \, \d t \right) \\
      & \leq \frac{1}{2 \pi} \int_0^{2 \pi} f(\gamma(e^{i t}))
      \, \d t.
    \end{split}
  \end{equation*}
  This shows that $f$ is \PSH. If equality holds, in the case where
  $f$ is strictly convex, $\gamma$ is essentially constant on the
  circle, and thus constant on $\Db$ by analytic continuation, so $f$
  is strictly \PSH.
\end{proof}

\begin{rmq}
  Theorem~\ref{thm:conv_psh} shows in particular that every continuous
  real linear form $f : E \ra \R$ is \PSH. Of course, this special
  case is a direct consequence of Cauchy's formula, which does not
  require Jensen's inequality.
\end{rmq}

\begin{thm}
  Let $E$ be a complete $\C$-normable space. Let $\sigma$ be a
  continuous anti-linear involution of $E$, and let $E_\R$ be the real
  subspace of $E$ fixed by $\sigma$, so that $E = E_\R \oplus i E_\R$.
  Denote by $\pi_\R : E \ra E_\R$ the projection along $i E_\R$ onto
  $E_\R$. Let $f_\R : E_\R \ra \R$ be a continuous function, and let
  $f = f_\R \circ \pi_\R$. Then $f_\R$ is convex if and only if $f$ is
  \PSH, and $f_\R$ is strictly convex if and only if $f$ is strictly
  \PSH.
\end{thm}

\begin{proof}
  Suppose that $f_\R$ is convex and let $\gamma : \Db \ra E$ be a
  holomorphic map. According to Cauchy's formula and Jensen's
  inequality,
  \begin{equation*}
    \begin{split}
      f(\gamma(0)) & = f \left( \frac{1}{2 \pi} \int_0^{2 \pi}
      \gamma(e^{i t}) \, \d t \right) \\
      & = f_\R \left( \frac{1}{2 \pi} \int_0^{2 \pi}
      \pi_\R(\gamma(e^{i t})) \, \d t \right) \\
      & \leq \frac{1}{2 \pi} \int_0^{2 \pi} f(\gamma(e^{i t}))
      \, \d t.
    \end{split}
  \end{equation*}
  This shows that $f$ is \PSH. If equality holds, in the case where
  $f_\R$ is strictly convex, $\pi_\R \circ \gamma$ is essentially
  constant on the circle, and thus constant by continuity, and this
  constant is $(\pi_\R \circ \gamma)(0)$. If $\lambda : E_\R \ra \R$
  is a continuous linear form, the same facts are true for $\lambda
  \circ \pi_\R \circ \gamma$. So this function, which is \PSH\ on
  $\D$, has a maximum point in $\D$; thus it is constant. We deduce
  that $\pi_\R \circ \gamma$ is itself constant on $\D$. Therefore, at
  every point of $\D$, the image of the differential of $\gamma$ is a
  complex subspace of $E$ contained in $i E_\R$; but such a subspace
  is necessarily zero, thus $\gamma$ is constant, so $f$ is strictly
  \PSH.

  Conversely, suppose that $f$ is \PSH\ and let $\gamma : \Sg \ra
  E_\R$ be an affine map. Let $\lambda : \Sg \ra \R$ be the linear map
  such that
  \begin{equation*}
    \lambda(1) = \frac{f_\R(\gamma(-1)) - f_\R(\gamma(1))}{2},
  \end{equation*}
  and let $g = f_\R \circ \gamma + \lambda$, so that
  \begin{align} \label{eq:g}
    g(0) & = f_\R(\gamma(0)), &
    g(-1) & = g(1) = \frac{f_\R(\gamma(-1)) + f_\R(\gamma(1))}{2}.
  \end{align}
  Let $U \subset \C$ be the open set defined by the condition $-1 <
  \real(z) < 1$, and let $\bar{U}$ be its closure. There exists a
  unique complex affine map $\tgamma : \bar{U} \ra E$ extending
  $\gamma$; moreover, $\pi_\R \circ \tgamma = \gamma \circ \real$, so
  $f \circ \tgamma = f_\R \circ \gamma \circ \real$. Thus the function
  $g \circ \real = f \circ \tgamma + \lambda \circ \real$ is \PSH\ on
  $U$. According to the maximum principle, thanks to the
  equalities~\eqref{eq:g},
  \begin{equation*}
    f_\R(\gamma(0))
    \leq \frac{f_\R(\gamma(-1)) + f_\R(\gamma(1))}{2}.
  \end{equation*}
  This shows that $f_\R$ is convex. If equality holds, then $g \circ
  \real$ is constant on $U$, so $f \circ \tgamma$ cannot be strictly
  \PSH; in the case where $f$ is strictly \PSH, we deduce that
  $\tgamma$ is constant, thus $\gamma$ is also constant, so $f_\R$ is
  strictly convex.
\end{proof}

\begin{cor}
  A continuous function $f : \C \ra \R$ which is constant on every
  vertical line is \PSH\ (resp. strictly \PSH) if and only if its
  restriction to $\R$ is convex (resp. strictly convex).
\end{cor}

\section{Spaces}

\subsection{Strictly convex spaces}

\begin{prp}
  Let $(E, \ndot)$ be a normed $\R$-vector space. The function $\ndot
  : E \ra \R$ is convex.
\end{prp}

\begin{proof}
  This is an immediate consequence of the triangle inequality.
\end{proof}

\begin{thm} \label{thm:conv_equiv}
  Let $(E, \ndot)$ be a normed $\R$-vector space. The following
  conditions are equivalent:
  \begin{enumerate}
  \item Every affine map $\gamma : \Sg \ra E$ whose image is contained
    in the unit sphere is constant.
  \item For all increasing, strictly convex maps $\psi : \R_+ \ra \R$,
    the composition $\psi \circ \ndot$ is strictly convex.
  \item There exists an increasing, strictly convex map $\psi : \R_+
    \ra \R$ such that the composition $\psi \circ \ndot$ is strictly
    convex.
  \end{enumerate}
\end{thm}

\begin{proof}
  Suppose condition~1, and let $\psi : \R_+ \ra \R$ be an increasing,
  strictly convex map and $\gamma : \Sg \ra E$ an affine map. By
  successively using the convexity of $\ndot$, the monotonicity of
  $\psi$ and the convexity of $\psi$, one obtains
  \begin{equation*}
    \begin{split}
      \psi(\norm{\gamma(0)})
      & \leq \psi \left( \frac{\norm{\gamma(-1)} +
        \norm{\gamma(1)}}{2} \right) \\
      & \leq \frac{\psi(\norm{\gamma(-1)}) +
        \psi(\norm{\gamma(1)})}{2}.
    \end{split}
  \end{equation*}
  If equality holds, by using the strict monotonicity and the strict
  convexity of $\psi$, one obtains $\norm{\gamma(0)} =
  \norm{\gamma(-1)} = \norm{\gamma(1)}$. Therefore, the function
  $\ndot \circ \gamma$, which is convex, has a maximum point in the
  interior of $\Sg$; thus it is constant. So the image of $\gamma$ is
  contained in a sphere, so $\gamma$ is constant. This proves
  condition~2.

  Condition~2 trivially implies condition~3; suppose that the latter
  is satisfied and let $\psi : \R_+ \ra \R$ be an increasing, strictly
  convex map such that $\psi \circ \ndot$ is strictly convex, and
  $\gamma : \Sg \ra E$ whose image is contained in the unit sphere of
  $\ndot$. The composition $\psi \circ \ndot \circ \gamma$ is
  constant, so $\gamma$ is constant by the strict convexity of $\psi
  \circ \ndot$, whence condition~1.
\end{proof}

\begin{dfn}
  If the conditions of Theorem~\ref{thm:conv_equiv} are satisfied,
  then $(E, \ndot)$ is said to be strictly convex.
\end{dfn}

\subsection{Strictly \psh\ spaces}

\begin{prp}
  Let $(E, \ndot)$ be a $\C$-Banach space. The function $\log \ndot :
  E \ra \Ri$ is \PSH.
\end{prp}

\begin{proof}
  Let $\gamma : \Db \ra \C$ be a holomorphic map, and for $r \in [0 ;
    1]$, let $N(r)$ be the number of zeros of $\gamma$ whose modulus
  is less than $r$. According to Jensen's formula,
  \begin{equation*}
    \log \abs{\gamma(0)} =
    \frac{1}{2 \pi} \int_0^{2 \pi} \log \abs{\gamma(e^{i t})}
    \, \d t - \int_0^1 \frac{N(r)}{r} \, \d r.
  \end{equation*}
  This proves the result when $E = \C$ and $\ndot = \abs{\cdot}$; now
  our aim is to reduce to this case. Let $\gamma : \Db \ra E$ be a
  holomorphic map. According to the Hahn-Banach theorem \cite[II,
    p.~67, cor.~1]{EVT}, there exists a continuous linear form
  $\lambda : E \ra \C$ of norm at most $1$ such that
  $\lambda(\gamma(0)) = \norm{\gamma(0)}$. We have
  \begin{equation*}
    \begin{split}
      \log \norm{\gamma(0)} & = \log \abs{\lambda(\gamma(0))} \\
      & \leq \frac{1}{2 \pi} \int_0^{2 \pi}
      \log \abs{\lambda(\gamma(e^{i t}))} \, \d t \\
      & \leq \frac{1}{2 \pi} \int_0^{2 \pi}
      \log \norm{\gamma(e^{i t})} \, \d t. \qedhere
    \end{split}
  \end{equation*}
\end{proof}

\begin{dfn}
  A map $\psi : \Ri \ra \Ri$ is said to be convex if it is continuous
  and its restriction to $\R$ is finite and convex. If in addition
  $\psi_{|\R}$ is strictly convex, then $\psi$ is said to be strictly
  convex.
\end{dfn}

\begin{rmq}
  Every convex (resp. strictly convex) map $\psi : \Ri \ra \Ri$ is
  increasing (resp. strictly increasing).
\end{rmq}

The next result is a variation of Jensen's inequality for functions
with values in $\Ri$.

\begin{lem} \label{lem:jensen_oo}
  Let $(S, \Sigma, \mu)$ be a measure space of total mass $1$ and $f :
  S \ra \Ri$ a measurable function bounded from above, and let $\psi :
  \Ri \ra \Ri$ be a convex map. We have
  \begin{equation*}
    \psi \left( \int_S f \d \mu \right)
    \leq \int_S \psi \circ f \d \mu.
  \end{equation*}
  Moreover, if $\psi$ is strictly convex and both sides are finite and
  equal, then $f$ is essentially constant.
\end{lem}

\begin{proof}
  Let $m = \int_S f \d \mu \in \Ri$. A slight extension of
  Lemma~\ref{lem:conv_aff} shows that there exists a continuous map
  $\alpha : \Ri \ra \Ri$ whose restriction to $\R$ is affine, such
  that $\alpha \leq \psi$ and $\alpha(m) = \psi(m)$. We have
  \begin{equation*}
    \begin{split}
      \psi \left( \int_S f \d \mu \right)
      & = \alpha \left( \int_S f \d \mu \right) \\
      & = \int_S \alpha \circ f \d \mu \\
      & \leq \int_S \psi \circ f \d \mu.
    \end{split}
  \end{equation*}
  If both sides are finite and equal, then $\alpha \circ f$ and $\psi
  \circ f$ coincide almost everywhere. Suppose that $\psi$ is strictly
  convex. If $m$ is finite, then $\alpha$ and $\psi$ only coincide at
  $m$ and possibly at $-\oo$; but in this case $f$ is finite almost
  everywhere, thus $f$ equals $m$ almost everywhere. If $m = -\oo$,
  then $\psi(-\oo) \in \R$ and $\alpha$ is constant, equal to
  $\psi(-\oo)$, so $\alpha$ and $\psi$ only coincide at $-\oo$, thus
  $f$ equals $-\oo$ almost everywhere.
\end{proof}

\begin{thm} \label{thm:psh_equiv}
  Let $(E, \ndot)$ be a $\C$-Banach space. The following conditions
  are equivalent:
  \begin{enumerate}
  \item Every holomorphic map $\gamma : \Db \ra E$ whose image is
    contained in the unit sphere is constant.
  \item For all strictly convex maps $\psi : \Ri \ra \Ri$, the
    composition $\psi \circ \log \ndot$ is strictly \PSH.
  \item There exists a strictly convex map $\psi : \Ri \ra \Ri$ such
    that the composition $\psi \circ \log \ndot$ is strictly \PSH.
  \end{enumerate}
\end{thm}

\begin{proof}
  Suppose condition~1, and let $\psi : \Ri \ra \Ri$ be a strictly
  convex map and $\gamma : \Db \ra E$ a holomorphic map. By
  successively using the fact that $\log \ndot$ is \PSH, the
  monotonicity of $\psi$ and Lemma~\ref{lem:jensen_oo}, one obtains
  \begin{equation*}
    \begin{split}
      \psi(\log \norm{\gamma(0)})
      & \leq \psi \left( \frac{1}{2 \pi} \int_0^{2 \pi}
      \log \norm{\gamma(e^{i t})} \, \d t \right) \\
      & \leq \frac{1}{2 \pi} \int_0^{2 \pi}
      \psi(\log \norm{\gamma(e^{i t})}) \, \d t.
    \end{split}
  \end{equation*}
  If equality holds, according to Lemma~\ref{lem:jensen_oo}, the
  function $\log \ndot \circ \gamma$ is essentially constant on the
  circle, thus constant by continuity; moreover, by using the strict
  monotonicity of $\psi$, we see that this constant is $\log
  \norm{\gamma(0)}$. Therefore, the function $\log \ndot \circ
  \gamma$, which is \PSH\ on $\D$, has a maximum point in $\D$; thus
  it is constant.  So the image of $\gamma$ is contained in a sphere,
  so $\gamma$ is constant. This proves condition~2.

  Condition~2 trivially implies condition~3; suppose that the latter
  is satisfied and let $\psi : \Ri \ra \Ri$ be a strictly convex map
  such that $\psi \circ \log \ndot$ is strictly \PSH, and $\gamma :
  \Db \ra E$ whose image is contained in the unit sphere of $\ndot$.
  The composition $\psi \circ \log \ndot \circ \gamma$ is constant, so
  $\gamma$ is constant since $\psi \circ \log \ndot$ is strictly \PSH,
  whence condition~1.
\end{proof}

\begin{dfn}
  If the conditions of Theorem~\ref{thm:psh_equiv} are satisfied, then
  $(E, \ndot)$ is said to be strictly \PSH.
\end{dfn}

\begin{rmq}
  Thorp and Whitley \cite[p.~641, th.~3.1]{TW67} proved that if
  condition~1 of Theorem~\ref{thm:psh_equiv} is satisfied for all
  \emph{affine} maps $\gamma : \Db \ra E$, then $(E, \ndot)$ is
  strictly \PSH.
\end{rmq}

\begin{prp}
  A $\C$-Banach space $(E, \ndot)$ which is strictly convex as an
  $\R$-Banach space is strictly \PSH.
\end{prp}

\begin{proof}
  Indeed, if there exists an increasing, strictly convex map $\psi :
  \R_+ \ra \R$ such that $\psi \circ \ndot$ is strictly convex, then
  $\psi \circ \ndot = \psi \circ \exp \circ \log \ndot$ is strictly
  \PSH\ according to Theorem~\ref{thm:conv_psh}; so we can conclude by
  using Theorem~\ref{thm:psh_equiv}, since $\psi \circ \exp : \Ri \ra
  \Ri$ is strictly convex.
\end{proof}

\subsection{\Pshe\ and \smmp}

\begin{dfn}
  A $\C$-Banach space $(E, \ndot)$ satisfies the \smmp\ if every
  holomorphic map $\eta$ from a connected holomorphic manifold $X$ to
  $E$ such that $\norm{\eta}$ has a local maximum is constant.
\end{dfn}

\begin{prp}
  A $\C$-Banach space $(E, \ndot)$ satisfies the \smmp\ if and only if
  it is strictly \PSH.
\end{prp}

\begin{proof}
  If $(E, \ndot)$ satisfies the \smmp, then clearly every holomorphic
  map $\gamma : \Db \ra E$ whose image is contained in the unit sphere
  is constant, so $(E, \ndot)$ is strictly \PSH.

  Conversely, suppose that $(E, \ndot)$ is strictly \PSH, let $X$ be a
  connected holomorphic manifold and $\eta : X \ra E$ a holomorphic
  map such that $\norm{\eta}$ has a local maximum at $x \in X$, and
  let us show that $\eta$ is constant.

  By analytic continuation, it is enough to show that $\eta$ is
  constant near $x$. Thus we can assume that $\norm{\eta}$ has a
  \emph{global} maximum at $x$, in which case this function is
  constant, since it is \PSH; moreover, we can assume that there
  exists a biholomorphism $\phi$ between a balanced open neighbourhood
  of the origin in a complete $\C$-normable space and $X$, satisfying
  $\phi(0) = x$. Now let $y \in X$. The holomorphic map $\gamma : \Db
  \ra E$ given by $\gamma(z) = \eta(\phi(z \cdot \phi^{-1}(y)))$ is
  constant, since its image is contained in a sphere; in particular,
  $\gamma(0) = \gamma(1)$, that is, $\eta(x) = \eta(y)$, which proves
  that $\eta$ is constant.
\end{proof}

\section{Direct integrals}

In this section, we let $\K$ denote $\R$ or $\C$.

\subsection{Measurability}

Let $(S, \Sigma, \mu)$ be a measure space and $\Ec = (E_s)_{s \in S}$
a family of $\K$-Banach spaces, and define
\begin{align*}
  \Phi_\Ec & = \{E_s \st s \in S\}, &
  \Pi_\Ec & = \prod_{s \in S} E_s, &
  \Omega_\Ec & = \coprod_{E \in \Phi_\Ec} E.
\end{align*}
Equip $\Phi_\Ec$ with the discrete $\sigma$-algebra, $\Pi_\Ec$ with
the product topology, and $\Omega_\Ec$ with the disjoint union
topology and with the corresponding Borel $\sigma$-algebra.

\begin{rms}
  $\Pi_\Ec$ is canonically identified with a closed subset of
  $(\Omega_\Ec)^S$. Also, $\Omega_\Ec$ is metrisable, as a disjoint
  union of metrisable spaces \cite[IX, p.~16]{TG_5-10}.
\end{rms}

\begin{dfn}
  If $\Ec$, seen as a map from $S$ to $\Phi_\Ec$, is measurable, then
  it is said to be a measurable family of $\K$-Banach spaces.
\end{dfn}

\begin{dfn}
  A section of $\Ec$ is an element of $\Pi_\Ec$. A measurable section
  of $\Ec$ is a section which, seen as a map from $S$ to $\Omega_\Ec$,
  is measurable.
\end{dfn}

\begin{dfn}
  Let $\sigma$ be a measurable section of $\Ec$, and let $\sigma(S)^*
  \subset \Omega_\Ec$ be the set of its non-zero values.
  \begin{itemize}
  \item If $\sigma(S)^*$ is finite, then $\sigma$ is said to be simple.
  \item If $\sigma(S)^*$ is separable, then $\sigma$ is said to be
    strongly measurable.
  \end{itemize}
\end{dfn}

\begin{prp}
  If the family $\Ec$ has a measurable section, then it is measurable.
\end{prp}

\begin{proof}
  Indeed, if $\sigma$ is a measurable section of $\Ec$, then for all
  subsets $X \subset \Phi_\Ec$, the inverse image of the open set
  $\coprod_{E \in X} E \subset \Omega_\Ec$ under $\sigma$ is
  $\Ec^{-1}(X)$; thus this set is measurable, so $\Ec$ is measurable.
\end{proof}

\begin{prp} \label{prp:mes_norm}
  For each $E \in \Phi_\Ec$, let $v_E$ be an element of $E$. For all
  measurable sections $\sigma$ of $\Ec$, the map $S \ra \R_+$ given by
  $s \mapsto \norm{\sigma(s) - v_{E_s}}_{E_s}$ is measurable.
\end{prp}

\begin{proof}
  For all $E \in \Phi_\Ec$, the map $\norm{\cdot - v_E}_E$ is
  continuous on $E$. The maps $\norm{\cdot - v_E}_E$ induce a
  continuous map $\Omega_\Ec \ra \R_+$, which is thus measurable. By
  composing this map with a measurable section $\sigma$ of $\Ec$, one
  obtains the map defined in the statement, which is therefore
  measurable.
\end{proof}

\begin{cor}
  For all measurable sections $\sigma$ of $\Ec$, the map $s \mapsto
  \norm{\sigma(s)}_{E_s}$ is measurable.
\end{cor}

\begin{prp}
  The set of measurable sections and the set of strongly measurable
  sections of $\Ec$ are sequentially closed in $\Pi_\Ec$.
\end{prp}

\begin{proof}
  Let $(\sigma_n)_{n \in \N}$ be a convergent sequence of measurable
  sections. Its limit $\sigma$ is measurable, as a pointwise limit of
  measurable functions with values in a metrisable space \cite[p.~245,
    prop.~8.1.10]{Cohn}.

  If in addition the $\sigma_n$ are strongly measurable, then for all
  $n \in \N$, $\sigma_n(S)^*$ is separable; thus the union $\bigcup_{n
    \in \N} \sigma_n(S)^*$ is also separable, as well as its closure;
  but the latter contains $\sigma(S)^*$, which is therefore separable.
  So $\sigma$ is strongly measurable.
\end{proof}

\begin{ntn}
  Given $E \in \Phi_\Ec$, $A \in \Sigma$ such that $\Ec$ is
  identically $E$ on $A$, and $v \in E$, we will denote by $\sigma_{A,
    v}$ the section of $\Ec$ which has value $v$ on $A$ and zero
  elsewhere.
\end{ntn}

\begin{thm}
  Suppose that $\Ec$ measurable.
  \begin{enumerate}
  \item The simple sections of $\Ec$ are the finite sums of sections
    of the form $\sigma_{A, v}$. In particular, the set of simple
    sections is a vector subspace of $\Pi_\Ec$.
  \item Every strongly measurable section $\sigma$ of $\Ec$ is the
    limit of a sequence $(\sigma_n)_{n \in \N}$ of simple sections
    such that for all $n \in \N$ and for all $s \in S$,
    $\norm{\sigma_n(s)}_{E_s} \leq \norm{\sigma(s)}_{E_s}$.

    In particular, the set of strongly measurable sections is the
    sequential closure of the set of simple sections; thus it is a
    vector subspace of $\Pi_\Ec$.
  \end{enumerate}
\end{thm}

\begin{proof}[Proof of assertion~1]
  For all simple sections $\sigma$,
  \begin{equation*}
    \sigma = \sum_{v \in \sigma(S)^*} \sigma_{A_v, v},
  \end{equation*}
  where $A_v = \sigma^{-1}(\{v\})$. Conversely, consider a section of
  the form
  \begin{equation*}
    \sigma = \sum_{i \in I} \sigma_{A_i, v_i},
  \end{equation*}
  where $I$ is a finite set. After rearranging, we can assume that the
  $A_i$ are non-empty and pairwise disjoint and that the $v_i$ are
  non-zero and pairwise distinct. Denote by $A$ the union of the
  $A_i$. Given any subset $Z \subset \Omega_\Ec$, we have
  \begin{equation*}
    \sigma^{-1}(Z) =
    \left( \bigcup_{v_i \in Z} A_i \right) \cup
    \big( \Ec^{-1}(\{E \in \Phi_\Ec \st 0_E \in Z\}) \, \bs \, A
    \big),
  \end{equation*}
  which shows that $\sigma$ is measurable. Moreover, $\sigma(S)^* =
  \{v_i \st i \in I\}$ is finite, so $\sigma$ is a simple section.
\end{proof}

\begin{proof}[Proof of assertion~2]
  Let $\sigma$ be a strongly measurable section of $\Ec$, let $D
  \subset \sigma(S)^*$ be a dense countable subset, and let $(v_n)_{n
    \in \N}$ be an enumeration of the set $\Q^\times D$ of non-zero
  rational multiples of elements of $D$. For $n \in \N$ and $s \in S$,
  define
  \begin{equation*}
    R_n(s) = \{k \in \lb 0 ; n \rb \st v_k \in E_s \text{ and }
    \norm{v_k}_{E_s} \leq \norm{\sigma(s)}_{E_s}\}.
  \end{equation*}
  If $R_n(s) = \varnothing$, define $\sigma_n(s) = 0_{E_s}$;
  otherwise, let $k$ be the smallest element of $R_n(s)$ minimising
  the quantity $\norm{\sigma(s) - v_k}_{E_s}$, and define $\sigma_n(s)
  = v_k$. Then $\norm{\sigma_n(s)}_{E_s} \leq \norm{\sigma(s)}_{E_s}$,
  and the sequence $(\sigma_n)_{n \in \N}$ converges to $\sigma$.

  Moreover, for all $n \in \N$ and for all $k \in \lb 0 ; n \rb$, the
  set $A_{n, k} = \sigma_n^{-1}(\{v_k\})$ is defined by the condition
  $v_k \in E_s$ and by inequalities involving the quantities
  $\norm{\sigma(s)}_{E_s}$ and $\norm{\sigma(s) - v_l}_{E_s}$, for $l
  \in \lb 0 ; n \rb$. Thus, thanks to Proposition~\ref{prp:mes_norm},
  it is measurable, so $\sigma_n = \sum_{k = 0}^n \sigma_{A_{n, k},
    v_k}$ is a simple section.
\end{proof}

\subsection{Integrability}

Henceforth, we assume that $\Ec$ is measurable, and we fix $p \in [1 ;
  \oo]$.

\begin{dfn}
  Let $\sigma$ be a measurable section of $\Ec$. We denote by
  $\norm{\sigma}_p$ the $p$-norm of the function $s \mapsto
  \norm{\sigma(s)}_{E_s}$. We say that $\sigma$ is $p$-integrable if
  $\norm{\sigma}_p < \oo$.
\end{dfn}

\begin{prp}
  The set $\Lc^p(\Ec)$ of strongly measurable $p$-integrable sections
  of $\Ec$ is a vector subspace of $\Pi_\Ec$, and $\ndot_p$ is a
  seminorm on $\Lc^p(\Ec)$, whose kernel is the space $\Nc(\Ec)$ of
  strongly measurable sections equal to zero almost everywhere.
\end{prp}

\begin{proof}
  The fact that the function $\ndot_p$ on the space of strongly
  measurable sections is positive, homogeneous and satisfies the
  triangle inequality comes from the corresponding properties for the
  norms $\ndot_{E_s}$ and for the $p$-norm of real measurable
  functions. This shows that $\Lc^p(\Ec)$ is a vector subspace of
  $\Pi_\Ec$ and that $\ndot_p$ is a seminorm on $\Lc^p(\Ec)$. The fact
  that the kernel of $\ndot_p$ is $\Nc(\Ec)$ again comes from the
  corresponding property for the $p$-norm of real measurable
  functions.
\end{proof}

\begin{dfn}
  The $L^p$ direct integral of the family $\Ec$, denoted by
  $L^p(\Ec)$, is the quotient $\Lc^p(\Ec) / \Nc(\Ec)$, equipped with
  the norm $\ndot_p$.
\end{dfn}

\begin{exm}
  If the family $\Ec$ is constant, equal to a Banach space $E$, then
  $L^p(\Ec)$ is the Bochner space $L^p(S, \Sigma, \mu ; E)$. In
  particular, if $E = \K$, then $L^p(\Ec)$ is the Lebesgue space
  $L^p(S, \Sigma, \mu ; \K)$.
\end{exm}

\begin{exm}
  Suppose that $\Sigma$ is discrete, that $\mu$ is the counting
  measure, and that the $E_s$ are pairwise distinct. Then $L^p(\Ec)$
  is the closed subspace $\ell^p_\den(\Ec)$ of $\ell^p(\Ec)$ whose
  elements are the sections with countable support. Note that
  $\ell^p_\den(\Ec) = \ell^p(\Ec)$ except when $p = \oo$ and the set
  of those $s \in S$ such that $E_s$ is non-zero is uncountable.
\end{exm}

\begin{exm}
  When all the $E_s$ are Hilbert spaces, $L^2(\Ec)$ is a Hilbert
  space. In particular, one recovers the notion of direct integral of
  Hilbert spaces, as it is defined in \cite[pp.~22-27]{N80}. Indeed,
  given a field $\Hc$ of Hilbert spaces along with a coherence
  $\alpha$, according to the terminology of that book, one obtains
  canonically a family $\Ec$ whose terms are equal to one of the
  spaces $\ell^2_n$, for $n \in \N$, or $\ell^2_\oo$. If $\Hc$ is
  Borel, then $\Ec$ is measurable, and $L^2(\mu ; \Hc, \alpha)$ is
  identified with $L^2(\Ec)$.
\end{exm}

The next result generalises the first two examples above.

\begin{thm}
  For $E \in \Phi_\Ec$, let $S_E = \Ec^{-1}(\{E\})$, and denote by
  $\Sigma_E$ and $\mu_E$ the restrictions of $\Sigma$ and $\mu$ to
  $S_E$. There exists a canonical isometric isomorphism
  \begin{equation*}
    L^p(\Ec) \simeq \ell^p_\den
    \big( (L^p(S_E, \Sigma_E, \mu_E ; E))_{E \in \Phi_\Ec} \big).
  \end{equation*}
\end{thm}

\begin{proof}
  First, for $\sigma \in \Pi_\Ec$ and $E \in \Phi_\Ec$, denote by
  $\sigma_E : S_E \ra E$ the restriction of $\sigma$ to $S_E$. The map
  $\sigma \mapsto (\sigma_E)_{E \in \Phi_\Ec}$ is an isomorphism
  \begin{equation*}
    \Pi_\Ec \simeq \prod_{E \in \Phi_\Ec} E^{S_E}.
  \end{equation*}

  Next, the image of the space of strongly measurable sections under
  this isomorphism is the space of families with countable support of
  strongly measurable sections.

  Indeed, if $\sigma$ is strongly measurable, then $(\sigma_E)_{E \in
    \Phi_\Ec}$ is a family of strongly measurable sections, which has
  countable support since $\sigma(S)^*$ is separable. Conversely, if
  $(\sigma_E)_{E \in \Phi_\Ec}$ is a family with countable support of
  strongly measurable sections, the measurability of $\Ec$ implies
  that for all $E \in \Phi_\Ec$, the section of $\Ec$ obtained by
  extending $\sigma_E$ by zero is strongly measurable. Therefore,
  $\sigma$ is strongly measurable, being the sum of a countable family
  of strongly measurable sections.

  Finally, if $\sigma$ is strongly measurable, then
  \begin{equation*}
    \norm{\sigma}_p = \norm[\big]{\big(
      \norm{\sigma_E}_p \big)_{E \in \Phi_\Ec}}_p.
  \end{equation*}
  Thus one obtains an isometric isomorphism
  \begin{equation*}
    \Lc^p(\Ec) \simeq \ell^p_\den \big( (\Lc^p(S_E, \Sigma_E, \mu_E ;
    E))_{E \in \Phi_\Ec} \big).
  \end{equation*}
  To conclude, it remains only to take the quotient in both sides by
  the kernel of the corresponding seminorm.
\end{proof}

\begin{dfn}
  The family $\Ec$ is said to be discrete if the image measure of
  $\mu$ under $\Ec$ is concentrated on the set of singletons of
  positive measure, that is, if $\mu(\Ec^{-1}(X)) = 0$, where $X = \{E
  \in \Phi_\Ec \st \mu(S_E) = 0\}$.
\end{dfn}

\begin{rmq}
  Suppose that $\mu$ is $\sigma$-finite and that $\Ec$ is not
  discrete. Then there exists a measurable subset of $\Ec^{-1}(X)$ of
  finite non-zero measure. The image measure under $\Ec$ of the
  restriction of $\mu$ to this subset is a finite non-zero measure on
  $\Phi_\Ec$, for which every singleton has measure zero. This implies
  that the cardinal of $\Phi_\Ec$ is greater than or equal to some
  inaccessible cardinal \cite[pp.~58-59]{F69}. Therefore, when $\mu$
  is $\sigma$-finite, the existence of a measurable, non-discrete
  family $\Ec$ cannot be proved in ZFC. Thus, in practice, all
  measurable families are discrete in this situation.
\end{rmq}

\begin{prp} \label{prp:inj}
  Suppose that $\mu$ is $\sigma$-finite and that $\Ec$ is discrete.
  Then for almost all $s \in S$, there exists an isometric embedding
  $E_s \ra L^p(\Ec)$.
\end{prp}

\begin{proof}
  One only has to show that every $E \in \Phi_\Ec$ such that $\mu(S_E)
  > 0$ can be imbedded into $L^p(\Ec)$. Now, if this condition is
  satisfied, since $\mu$ is $\sigma$-finite, there exists a measurable
  subset $A \subset S_E$ of finite non-zero measure, and the map which
  associates to each $v \in E$ the class of
  $\mu(A)^{-\frac{1}{p}}\sigma_{A, v}$ in $L^p(\Ec)$ is an isometric
  embedding.
\end{proof}

\subsection{Completeness of direct integrals}

\begin{dfn}
  Let $X$ be a separable topological space, and let $\sigma : X \ra
  \Lc^p(\Ec)$ be a map such that $\sigma(\cdot)(s) : X \ra E_s$ is
  continuous for all $s \in S$. We call $\sigma$ a parametric section,
  and we denote by $\norm{\sigma}_p$ the $p$-norm of the function $s
  \mapsto \sup_{x \in X} \norm{\sigma(x)(s)}_{E_s}$.
\end{dfn}

\begin{rmq}
  This function is measurable, because if $(x_n)_{n \in \N}$ is a
  dense sequence of points of $X$, then $\sup_{x \in X}
  \norm{\sigma(x)(s)}_{E_s} = \sup_{n \in \N}
  \norm{\sigma(x_n)(s)}_{E_s}$.
\end{rmq}

\begin{thm} \label{thm:banach}
  Let $(\sigma_n)_{n \in \N}$ be a sequence of parametric sections
  such that $\sum_{n = 0}^\oo \norm{\sigma_n}_p < \oo$. There exists a
  parametric section $\sigma : X \ra \Lc^p(\Ec)$ such that
  \begin{equation*}
    \lim_{N \ra \oo} \norm[\bigg]{\sigma - \sum_{n = 0}^N \sigma_n}_p
    = 0.
  \end{equation*}
  Moreover, for almost all $s \in S$, the series of functions $\sum_{n
    = 0}^\oo \sigma_n(\cdot)(s) : X \ra E_s$ converges normally to
  $\sigma(\cdot)(s)$.
\end{thm}

\begin{proof}
  Let $f : S \ra [0 ; \oo]$ be the measurable function defined by
  \begin{equation*}
    f(s) =
    \sum_{n = 0}^\oo \sup_{x \in X} \norm{\sigma_n(x)(s)}_{E_s}.
  \end{equation*}
  We have $\norm{f}_p \leq \sum_{n = 0}^\oo \norm{\sigma_n}_p < \oo$,
  so $f$ is finite almost everywhere. For $x \in X$ and $s \in S$,
  define
  \begin{equation*}
    \sigma(x)(s) = \left\{
    \begin{aligned}
      & \sum_{n = 0}^\oo \sigma_n(x)(s)
      & \quad & \text{if $f(s) < \oo$}, \\
      & 0 & \quad & \text{otherwise}.
    \end{aligned}
    \right.
  \end{equation*}
  For all $s \in S$ such that $f(s) < \oo$, we have
  \begin{equation*}
    \sup_{x \in X} \norm[\bigg]{\sigma(x)(s) - \sum_{n = 0}^N
      \sigma_n(x)(s)}_{E_s} \leq
    \sum_{n = N + 1}^\oo \sup_{x \in X} \norm{\sigma_n(x)(s)}_{E_s}.
  \end{equation*}
  Therefore,
  \begin{equation*}
    \begin{split}
      \norm[\bigg]{\sigma - \sum_{n = 0}^N \sigma_n}_p
      & \leq \sum_{n = N + 1}^\oo \norm{\sigma_n}_p \\
      & \xrightarrow[N \ra \oo]{} 0.
    \end{split}
  \end{equation*}
  Finally, $\sigma$ is a parametric section having all the properties
  claimed.
\end{proof}

\begin{cor}
  $L^p(\Ec)$ is a $\K$-Banach space.
\end{cor}

\begin{proof}
  Indeed, Theorem~\ref{thm:banach} applied to a one-point space $X$
  shows that every absolutely convergent series of elements of
  $L^p(\Ec)$ is convergent, from which the result follows according to
  \cite[IX, p.~37, cor.~2]{TG_5-10}.
\end{proof}

\subsection{Strict convexity of direct integrals}

In this subsection, we assume that $\K = \R$.

\begin{lem} \label{lem:conv_tgamma}
  For every affine map $\gamma : \Sg \ra L^p(\Ec)$, there exists a
  parametric section $\tgamma : \Sg \ra \Lc^p(\Ec)$ lifting $\gamma$
  such that for all $s \in S$, the map $\tgamma(\cdot)(s) : \Sg \ra
  E_s$ is affine.
\end{lem}

\begin{proof}
  There exist elements $\gamma_0$ and $\gamma_1$ of $L^p(\Ec)$ such
  that for all $t \in \Sg$,
  \begin{equation*}
    \gamma(t) = \gamma_0 + t \gamma_1.
  \end{equation*}
  Let $\tgamma_0$ and $\tgamma_1$ be elements of $\Lc^p(\Ec)$ lifting
  $\gamma_0$ and $\gamma_1$, and for $t \in \Sg$, define
  \begin{equation*}
    \tgamma(t) = \tgamma_0 + t \tgamma_1.
  \end{equation*}
  Then $\tgamma$ is a lifting of $\gamma$ satisfying the required
  conditions.
\end{proof}

\begin{thm} \label{thm:conv_int}
  Suppose that $1 < p < \oo$.
  \begin{enumerate}
  \item If $E_s$ is strictly convex for almost all $s \in S$, then
    $L^p(\Ec)$ is strictly convex.
  \item Suppose that $\mu$ is $\sigma$-finite and that $\Ec$ is
    discrete. If $L^p(\Ec)$ is strictly convex, then $E_s$ is strictly
    convex for almost all $s \in S$.
  \end{enumerate}
\end{thm}

\begin{proof}
  Suppose that $E_s$ is strictly convex for almost all $s$, and let us
  show that $L^p(\Ec)$ is strictly convex. According to
  Theorem~\ref{thm:conv_equiv}, this amounts to assuming that
  $\ndot_{E_s}^p$ is strictly convex for almost all $s$, and we have
  to show that $\ndot_p^p$ is strictly convex.

  Let $\gamma : \Sg \ra L^p(\Ec)$ be an affine map, and let $\tgamma$
  be a lifting of $\gamma$ satisfying the statement of
  Lemma~\ref{lem:conv_tgamma}. We have
  \begin{equation*}
    \begin{split}
      \norm{\gamma(0)}_p^p
      & = \int_S \norm{\tgamma(0)(s)}_{E_s}^p \d \mu(s) \\
      & \leq \int_S \frac{\norm{\tgamma(-1)(s)}_{E_s}^p +
        \norm{\tgamma(1)(s)}_{E_s}^p}{2} \, \d \mu(s) \\
      & = \frac{\norm{\gamma(-1)}_p^p + \norm{\gamma(1)}_p^p}{2}.
    \end{split}
  \end{equation*}
  If equality holds, then for almost all $s$,
  \begin{equation*}
    \norm{\tgamma(0)(s)}_{E_s}^p =
    \frac{\norm{\tgamma(-1)(s)}_{E_s}^p +
      \norm{\tgamma(1)(s)}_{E_s}^p}{2},
  \end{equation*}
  thus $\tgamma(\cdot)(s)$ is constant for almost all $s$, so $\gamma$
  is constant. This proves assertion~1. Assertion~2 is a consequence
  of Proposition~\ref{prp:inj}.
\end{proof}

\subsection{Strict \pshe\ of direct integrals}

In this subsection, we assume that $\K = \C$.

\begin{lem} \label{lem:psh_tgamma}
  For every holomorphic map $\gamma : \Db \ra L^p(\Ec)$, there exists
  a parametric section $\tgamma : \Db \ra \Lc^p(\Ec)$ lifting $\gamma$
  such that for almost all $s \in S$, the map $\tgamma(\cdot)(s) : \Db
  \ra E_s$ is holomorphic.
\end{lem}

\begin{proof}
  There exists a sequence $(\gamma_n)_{n \in \N}$ of elements of
  $L^p(\Ec)$ and a real number $r > 1$ such that $\sum_{n = 0}^\oo r^n
  \norm{\gamma_n}_p < \oo$ and for all $z \in \Db$,
  \begin{equation*}
    \gamma(z) = \sum_{n = 0}^\oo z^n \gamma_n.
  \end{equation*}
  Let $(\tgamma_n)_{n \in \N}$ be a sequence of elements of
  $\Lc^p(\Ec)$ lifting $(\gamma_n)_{n \in \N}$. According to
  Theorem~\ref{thm:banach}, applied to the topological space $r \Db$
  and to the parametric sections $z \mapsto z^n \tgamma_n$, there
  exists a parametric section $\tgamma : \Db \ra \Lc^p(\Ec)$ lifting
  $\gamma$ such that for almost all $s \in S$, $\sum_{n = 0}^\oo r^n
  \norm{\tgamma_n(s)}_{E_s} < \oo$, and for all $z \in \Db$,
  \begin{equation*}
    \tgamma(z)(s) = \sum_{n = 0}^\oo z^n \tgamma_n(s).
  \end{equation*}
  Then $\tgamma$ is a lifting of $\gamma$ satisfying the required
  conditions.
\end{proof}

\begin{lem}[{\cite[p.~122, lemma 9.2]{M52}}] \label{lem:mes_prod}
  Let $K$ be a compact metric space and $Y$ a metrisable topological
  space, equipped with their Borel $\sigma$-algebras $\Bc(K)$ and
  $\Bc(Y)$, and let $f : K \times S \ra Y$. Suppose that $f(x, \cdot)
  : S \ra Y$ is measurable for all $x \in K$ and that $f(\cdot, s) : K
  \ra Y$ is continuous for all $s \in S$. Then $f$ is measurable with
  respect to the product $\sigma$-algebra $\Bc(K) \times \Sigma$.
\end{lem}

\begin{proof}
  For $n \geq 1$, let $K = \coprod_{i \in I_n} K_{n, i}$ be a finite
  partition of $K$ into Borel subsets of diameter less than
  $\frac{1}{n}$, and for $i \in I_n$, let $x_{n, i} \in K_{n, i}$.
  Let $f_n : K \times S \ra Y$ be defined by $f_n(x, s) = f(x_{n, i},
  s)$ for $(x, s) \in K_{n, i} \times S$. Then $f_n$ is measurable
  with respect to $\Bc(K) \times \Sigma$, and the sequence $(f_n)_{n
    \geq 1}$ converges pointwise to $f$, which implies the desired
  result according to \cite[p.~245, prop.~8.1.10]{Cohn}.
\end{proof}

\begin{thm} \label{thm:psh_int}
  Suppose that $\mu$ is $\sigma$-finite and that $1 \leq p < \oo$.
  \begin{enumerate}
  \item If $E_s$ is strictly \PSH\ for almost all $s \in S$, then
    $L^p(\Ec)$ is strictly \PSH.
  \item Suppose that $\Ec$ is discrete. If $L^p(\Ec)$ is strictly
    \PSH, then $E_s$ is strictly \PSH\ for almost all $s \in S$.
  \end{enumerate}
\end{thm}

\begin{proof}
  Suppose that $E_s$ is strictly \PSH\ for almost all $s$, and let us
  show that $L^p(\Ec)$ is strictly \PSH. According to
  Theorem~\ref{thm:psh_equiv}, this amounts to assuming that
  $\ndot_{E_s}^p$ is strictly \PSH\ for almost all $s$, and we have to
  show that $\ndot_p^p$ is strictly \PSH.

  Let $\gamma : \Db \ra L^p(\Ec)$ be a holomorphic map, and let
  $\tgamma$ be a lifting of $\gamma$ satisfying the statement of
  Lemma~\ref{lem:psh_tgamma}. We have
  \begin{equation*}
    \begin{split}
      \norm{\gamma(0)}_p^p
      & = \int_S \norm{\tgamma(0)(s)}_{E_s}^p \d \mu(s) \\
      & \leq \int_S \frac{1}{2 \pi} \int_0^{2 \pi}
      \norm{\tgamma(e^{i t})(s)}_{E_s}^p \d t \, \d \mu(s) \\
      & = \frac{1}{2 \pi} \int_0^{2 \pi}
      \norm{\gamma(e^{i t})}_p^p \, \d t,
    \end{split}
  \end{equation*}
  thanks to Lemma~\ref{lem:mes_prod}, applied to the function $(t, s)
  \ra \norm{\tgamma(e^{i t})(s)}_{E_s}^p$, and to Fubini's theorem.
  If equality holds, then for almost all $s$,
  \begin{equation*}
    \norm{\tgamma(0)(s)}_{E_s}^p =
    \frac{1}{2 \pi} \int_0^{2 \pi}
    \norm{\tgamma(e^{i t})(s)}_{E_s}^p \d t,
  \end{equation*}
  thus $\tgamma(\cdot)(s)$ is constant for almost all $s$, so $\gamma$
  is constant. This proves assertion~1. Assertion~2 is a consequence
  of Proposition~\ref{prp:inj}.
\end{proof}

\appendix

\section{Day's method}

The purpose of this appendix is to prove Theorems \ref{thm:conv_int}
and \ref{thm:psh_int} again, under the hypothesis $1 < p < \oo$, by
Day's method \cite[p.~520, th.~6]{D55}. In both cases, we will not
come back to assertion~2, which is a consequence of
Proposition~\ref{prp:inj}.

\begin{proof}[Proof of Theorem \ref{thm:conv_int}]
  Let us begin by observing that a normed $\R$-vector space $(E,
  \ndot)$ is strictly convex if and only if for all affine maps
  $\gamma : \Sg \ra E$, the equalities $\norm{\gamma(0)} =
  \norm{\gamma(-1)} = \norm{\gamma(1)}$ imply $\gamma(-1) =
  \gamma(1)$.

  Now let $\gamma : \Sg \ra L^p(\Ec)$ be an affine map, let $\tgamma$
  be a lifting of $\gamma$ satisfying the statement of
  Lemma~\ref{lem:conv_tgamma}, let $\trho : \Sg \ra \Lc^p(S, \Sigma,
  \mu ; \R)$ be the map defined by $\trho(t)(s) =
  \norm{\tgamma(t)(s)}_{E_s}$, and let $\rho : \Sg \ra L^p(S, \Sigma,
  \mu ; \R)$ be the composition of $\trho$ and the quotient map. For
  all $s \in S$, since $\tgamma(\cdot)(s)$ is affine and $\ndot_{E_s}$
  is convex,
  \begin{equation} \label{eq:conv_trho}
    \trho(0)(s) \leq \frac{\trho(-1)(s) + \trho(1)(s)}{2}.
  \end{equation}
  Therefore,
  \begin{equation*}
    \begin{split}
      \norm{\gamma(0)}_p & = \norm{\rho(0)}_p \\
      & \leq \norm[\bigg]{\frac{\rho(-1) + \rho(1)}{2}}_p \\
      & \leq \frac{\norm{\rho(-1)}_p + \norm{\rho(1)}_p}{2} \\
      & = \frac{\norm{\gamma(-1)}_p + \norm{\gamma(1)}_p}{2}.
    \end{split}
  \end{equation*}
  If $\norm{\gamma(0)}_p = \norm{\gamma(-1)}_p = \norm{\gamma(1)}_p$,
  then we have equality, so the inequality~\eqref{eq:conv_trho} is an
  equality for almost all $s$; moreover,
  \begin{equation} \label{eq:conv_rho}
    \norm[\bigg]{\frac{\rho(-1) + \rho(1)}{2}}_p = \norm{\rho(-1)}_p =
    \norm{\rho(1)}_p.
  \end{equation}
  Since $L^p(S, \Sigma, \mu ; \R)$ is strictly convex, the
  equalities~\eqref{eq:conv_rho} imply $\rho(-1) =
  \rho(1)$. Therefore, for almost all $s$, $\trho(0)(s) = \trho(-1)(s)
  = \trho(1)(s)$, that is, $\norm{\tgamma(0)(s)}_{E_s} =
  \norm{\tgamma(-1)(s)}_{E_s} = \norm{\tgamma(1)(s)}_{E_s}$. If $E_s$
  is strictly convex for almost all $s$, this shows that $\gamma(-1) =
  \gamma(1)$, so $L^p(\Ec)$ is strictly convex.
\end{proof}

\begin{lem} \label{lem:jensen_norm}
  Let $(E, \ndot)$ be a strictly convex $\R$-Banach space, $(S,
  \Sigma, \mu)$ a measure space of total mass $1$, and $\eta : S \ra
  E$ an integrable map. We have
  \begin{equation*}
    \norm[\bigg]{\int_S \eta \, \d \mu}
    \leq \int_S \norm{\eta} \, \d \mu.
  \end{equation*}
  If equality holds and $\norm{\eta}$ is constant, then $\eta$ is
  essentially constant.
\end{lem}

\begin{proof}
  The inequality is clear. Suppose that equality holds and that
  $\norm{\eta}$ is constant, and let $\psi : \R_+ \ra \R$ be an
  increasing, strictly convex map. We have
  \begin{equation*}
    \begin{split}
      \psi \left( \norm[\bigg]{\int_S \eta \, \d \mu} \right)
      & = \psi \left( \int_S \norm{\eta} \, \d \mu \right) \\
      & = \int_S \psi(\norm{\eta}) \, \d \mu.
    \end{split}
  \end{equation*}
  Lemma~\ref{lem:jensen} then shows that $\eta$ is essentially
  constant, since $\psi \circ \ndot$ is strictly convex.
\end{proof}

\begin{proof}[Proof of Theorem \ref{thm:psh_int}]
  Let us begin by observing that a $\C$-Banach space $(E, \ndot)$ is
  strictly \PSH\ if and only if for all holomorphic maps $\gamma : \Db
  \ra E$, the equalities $\norm{\gamma(0)} = \norm{\gamma(z)}$, for
  $\abs{z} = 1$, imply that $\gamma$ is constant on the circle.

  Now let $\gamma : \Db \ra L^p(\Ec)$ be a holomorphic map, let
  $\tgamma$ be a lifting of $\gamma$ satisfying the statement of
  Lemma~\ref{lem:psh_tgamma}, let $\trho : \Db \ra \Lc^p(S, \Sigma,
  \mu ; \R)$ be the map defined by $\trho(z)(s) =
  \norm{\tgamma(z)(s)}_{E_s}$, and let $\rho : \Db \ra L^p(S, \Sigma,
  \mu ; \R)$ be the composition of $\trho$ and the quotient map. For
  almost all $s \in S$, since $\tgamma(\cdot)(s)$ is holomorphic and
  $\ndot_{E_s}$ is \PSH,
  \begin{equation} \label{eq:psh_trho}
    \trho(0)(s) \leq
    \frac{1}{2 \pi} \int_0^{2 \pi} \trho(e^{i t})(s) \, \d t.
  \end{equation}
  Therefore, thanks to Lemma~\ref{lem:mes_prod}, applied to the
  function $(t, s) \mapsto \trho(e^{i t})(s)$, and to the result
  \cite[p.~26, prop.~1.2.25]{HNVW16} about pointwise calculation of
  integrals with values in an $L^p$ space,
  \begin{equation*}
    \begin{split}
      \norm{\gamma(0)}_p & = \norm{\rho(0)}_p \\
      & \leq \norm[\bigg]{\frac{1}{2 \pi}
        \int_0^{2 \pi} \rho(e^{i t}) \, \d t}_p \\
      & \leq \frac{1}{2 \pi}
      \int_0^{2 \pi} \norm{\rho(e^{i t})}_p \d t \\
      & = \frac{1}{2 \pi}
      \int_0^{2 \pi} \norm{\gamma(e^{i t})}_p \d t.
    \end{split}
  \end{equation*}
  If $\norm{\gamma(0)}_p = \norm{\gamma(z)}_p$ for $\abs{z} = 1$, then
  we have equality, so the inequality~\eqref{eq:psh_trho} is an
  equality for almost all $s$; moreover, for $\abs{z} = 1$,
  \begin{equation} \label{eq:psh_rho}
    \norm[\bigg]{\frac{1}{2 \pi}
      \int_0^{2 \pi} \rho(e^{i t}) \, \d t}_p = \norm{\rho(z)}_p.
  \end{equation}
  Since $L^p(S, \Sigma, \mu ; \R)$ is strictly convex, the
  equalities~\eqref{eq:psh_rho} and Lemma~\ref{lem:jensen_norm} imply
  that $\rho$ is essentially constant on the circle. Thus there exists
  a dense sequence $(z_n)_{n \in \N}$ of points of the circle such
  that for all $n \in \N$, for almost all $s \in S$, $\trho(z_n)(s) =
  \trho(z_0)(s)$. We deduce that for almost all $s$, for all $n \in
  \N$, $\trho(z_n)(s) = \trho(z_0)(s)$; so, by continuity, for almost
  all $s$, for $\abs{z} = 1$, $\trho(0)(s) = \trho(z)(s)$, that is,
  $\norm{\tgamma(0)(s)}_{E_s} = \norm{\tgamma(z)(s)}_{E_s}$. If $E_s$
  is strictly \PSH\ for almost all $s$, this shows that $\gamma$ is
  constant on the circle, so $L^p(\Ec)$ is strictly \PSH.
\end{proof}

\begin{rmq}
  In the proof of Theorem~\ref{thm:conv_int}, the strict convexity of
  $L^p(S, \Sigma, \mu ; \R)$ and the equalities~\eqref{eq:conv_rho}
  enable one to conclude that $\rho(-1) = \rho(1)$; the reason is that
  there exists an affine map $\Sg \ra L^p(S, \Sigma, \mu ; \R)$
  agreeing with $\rho$ at $-1$ and $1$. In the proof of
  Theorem~\ref{thm:psh_int}, if there existed a holomorphic map $\Db
  \ra L^p(S, \Sigma, \mu ; \C)$ agreeing with $\rho$ on the circle,
  one could use the strict \pshe\ of this space and the
  equalities~\eqref{eq:psh_rho} to conclude that $\rho$ is constant on
  the circle, in which case the reasoning would remain valid for $p =
  1$.  Unfortunately, this is not the case \emph{a priori}, which
  suggests that Day's method cannot be used to prove
  Theorem~\ref{thm:psh_int} in full generality.
\end{rmq}

\bibliography{biblio}

\signature{Anne-Edgar \textsc{Wilke}, Univ. Bordeaux, CNRS, INRIA,
  Bordeaux INP, IMB, UMR 5251, F-33400 Talence, France}

\end{document}